\documentclass[12pt]{article}
\usepackage{mathrsfs}
\usepackage{amssymb,amsfonts,amsmath,amsthm,cite}
\usepackage{epsfig}
\newtheorem{thm}{Theorem}[section]

\newtheorem{lem}[thm]{Lemma}

\makeatletter \@addtoreset{equation}{section}

\def\qed{\hfill \rule{4pt}{7pt}}

\begin{document}

{\Large
\begin{center}
Interlacing Log-concavity of the Boros-Moll Polynomials

\end{center}
}

\begin{center}
William Y. C. Chen$^1$,  Larry X. W. Wang$^2$
and Ernest X. W. Xia$^3$ \\

Center for Combinatorics, LPMC-TJKLC\\
Nankai University\\
 Tianjin 300071, P. R. China

$^1$chen@nankai.edu.cn, $^2$wxw@cfc.nankai.edu.cn,\\
 $^3$xxwrml@mail.nankai.edu.cn

\end{center}


\noindent {\bf Abstract.} We introduce the notion of interlacing log-concavity of a polynomial sequence $\{P_m(x)\}_{m\geq 0}$, where $P_m(x)$ is a polynomial of degree $m$ with positive coefficients $a_{i}(m)$. This sequence of polynomials is said to be interlacing log-concave if the ratios of consecutive coefficients of $P_m(x)$  interlace the ratios of consecutive coefficients of $P_{m+1}(x)$ for any $m\geq 0$.
Interlacing log-concavity is stronger than the log-concavity. We show
that the Boros-Moll polynomials are  interlacing
log-concave. Furthermore we give a
sufficient condition for
interlacing log-concavity which implies that some classical combinatorial polynomials are interlacing log-concave.



\section{Introduction}

In this paper, we introduce the notion of interlacing log-concavity of a polynomial sequence $\{P_m(x)\}_{m\geq 0}$, which is stronger than the log-concavity of the polynomials $P_m(x)$.
We shall show that the Boros-Moll polynomials are
interlacing log-concave.

For a sequence polynomials $\{P_m(x)\}$, 
let \[ P_m(x)=\sum \limits_{i=0}^{m}a_i(m)x^m,\]
and let $r_i(m) = a_i(m) / a_{i+1}(m)$. We say that the polynomials $P_m(x)$  are   interlacing log-concave  if the ratios $r_i(m)$ interlace the ratios $r_i(m+1)$, that is,
\begin{align}
r_0(m+1) \leq r_0(m) \leq r_1(m+1) \leq r_1(m) \leq \cdots \leq r_{m-1}(m+1)\leq r_{m-1}(m)\leq r_m(m+1).
\end{align}

Recall that a sequence $\{a_i\}_{0 \leq i \leq  m}$
 of positive numbers is said to be log-concave if
\[
\frac{a_0}{a_1} \leq \frac{a_1}{a_2} \leq \cdots  \leq
\frac{a_{m-1}}{a_m}.
\]
It is clear that the interlacing log-concavity implies the log-concavity.

 For the background on the Boros-Moll polynomials;
see \cite{Amdeberhan2007,George1999-1,George1999-2,
 George1999-3,George2001
 ,George2004,Moll2002}.
From now on, we shall use $P_m(a)$ to denote
the Boros-Moll polynomial given by
\begin{align}
P_m(x)=\sum_{j,k}{2m+1 \choose 2j}{m-j \choose k}{2k+2j \choose
k+j}\frac{(x+1)^j(x-1)^k}{2^{3(k+j)}}.
\end{align}
Boros and Moll \cite{George1999-1} derived the following formula for the coefficient $d_i(m)$ of $x^i$ in $P_m(x)$,
\begin{align}\label{Defi}
d_i(m)=2^{-2m}\sum_{k=i}^m2^k{2m-2k \choose m-k}{m+k \choose
k}{k\choose i}.
\end{align}

  Boros and Moll \cite{George1999-2}
 proved that
the sequence $\{d_i(m)\}_{0 \leq i \leq m}$ is unimodal and the
maximum element appears in the middle. In other words,
\begin{align}\label{M}
d_0(m)<d_1(m)< \cdots <d_{\left[\frac{m}{2}\right]}(m)>
d_{\left[\frac{m}{2}\right]-1}(m)>\cdots >d_m(m).
\end{align}
   Moll \cite{Moll2002}  conjectured 
   $P_m(x)$ is log-concave for any $m$. Kauers and Paule \cite{Kauers2006} confirmed this
conjecture based on recurrence relations found by a  computer
algebra approach. Chen and Xia \cite{Chen2008}  showed that the
sequence
  $\{d_i(m)\}_{0 \leq i \leq m}$ satisfies the ratio monotone
  property which implies the log-concavity
  and the spiral property. Chen and Gu showed that for any $m$,  $P_m(x)$ is reverse ultra log-concave \cite{CG09}.

The main result of this paper is to show that the Boros-Moll polynomials
are  interlacing log-concave. We also give a sufficient condition for the
interlacing log-concavity from which we see that several classical
combinatorial polynomials are interlacing log-concave.

\section{The interlacing   log-concavity of
$d_i(m)$} \label{section3}

In this section, we show that for $m\geq 2$, the
 the Boros-Moll  polynomials $P_m(x)$ are
interlacing log-concave. More precisely, we have

\begin{thm}\label{Theorem1}
For $m\geq 2$ and $0\leq i \leq m$, we have
\begin{equation}\label{cd1}
d_i(m)d_{i+1}(m+1)>d_{i+1}(m)d_i(m+1)
\end{equation}
and
\begin{equation}\label{cd2}
d_i(m)d_i(m+1)>d_{i-1}(m)d_{i+1}(m+1).
\end{equation}

\end{thm}

The proof relies on
the following recurrence relations derived by Kauers and Paule
\cite{Kauers2006}. In fact, they found four recurrence relations
 for the Boros-Moll sequence  $\{d_i(m)\}_{0 \leq i \leq m}$:
 \allowdisplaybreaks
\begin{align}
d_i(m+1)=&\frac{m+i}{m+1}d_{i-1}(m) +\frac{(4m+2i+3)}{2(m+1)}d_i(m),
 \ \ \ \ 0 \leq i \leq m+1, \label{recu1}\\[6pt]
 d_{i}(m+1)=&\frac{(4m-2i+3)(m+i+1)}{2(m+1)(m+1-i)}d_i(m)
 \nonumber \\[6pt]
 & \quad -\frac{i(i+1)}{(m+1)(m+1-i)}d_{i+1}(m),
 \qquad \ \ \ 0 \leq i \leq
 m, \label{recu2}\\[6pt]
 d_i(m+2)=&\frac{-4i^2+8m^2+24m+19}{2(m+2-i)
 (m+2)}d_i(m+1) \nonumber \\[6pt]
 & \quad -\frac{(m+i+1)(4m+3)(4m+5)}
 {4(m+2-i)(m+1)(m+2)}d_i(m), \qquad \ 0 \leq
  i \leq m+1,\label{recu3}
\end{align}
and for $0 \leq i \leq m+1$,
\begin{align}\label{recu4}
(m+2-i)(m+i-1)d_{i-2}(m)-(i-1)(2m+1)d_{i-1}(m)+i(i-1)d_i(m)=0.
\end{align}
Note that Moll \cite{Moll2007} also has independently derived the recurrence relation (\ref{recu4}) from which  the other three   relations can be deduced.

To prove (\ref{cd1}), we  give the following lemma.

\begin{lem}\label{strlog}Let $m \geq 2$ be an integer.
For $0\leq i\leq m-2$, we have
\begin{equation}\label{W-result}
\frac{d_i(m)} {d_{i+1}(m)} <
\frac{(4m+2i+3)d_{i+1}(m)}{(4m+2i+7)d_{i+2}(m)}
.
\end{equation}
\end{lem}

\proof We proceed by induction on $m$.   It is easy to
check that the theorem is valid for $m=2$. Assume that the result is true for  $
n$, that is, for $0\leq i\leq n-2$,
\begin{align}\label{W-assume}
\frac{d_i(n)} {d_{i+1}(n)}< \frac{(4n+2i+3)d_{i+1}(n)}{(4n+2i+7)d_{i+2}(n)}.
\end{align}
We aim to show that \eqref{W-result} holds for $n+1$, that is, for $0 \leq i \leq n-1$,
\begin{equation}\label{W-result-1}
\frac{d_i(n+1)} {d_{i+1}(n+1)} <
\frac{(4n+2i+7)d_{i+1}(n+1)}{(4n+2i+11)
d_{i+2}(n+1)}
 .
\end{equation}
From the recurrence relation \eqref{recu1},
  we can verify that for $0 \leq i \leq n-1$,
\begin{align*}
& (2i+4n+7)d_{i+1}^2(n+1)
-(2i+4n+11)d_i(n+1)d_{i+2}(n+1)\\[6pt]
&\quad =(2i+4n+7)\left(\frac{i+n+1}{n+1}d_i(n)+\frac{2i+4n+5}
{2(n+1)}d_{i+1}(n)\right)^2\\[6pt]
& \qquad -(2i+4n+11)\left(\frac{i+n+2}{n+1}d_{i+1}(n)
+\frac{2i+4n+7}{2(n+1)}d_{i+2}(n)
\right)\\[6pt]
& \qquad \quad \times \left(\frac{n+i}{n+1}d_{i-1}(n)
+\frac{2i+4n+3}{2(n+1)}d_i(n)\right)\\[6pt]
&\quad =\frac{A_1(n,i)+A_2(n,i)+A_3(n,i)}{{4(n+1)^2}},
\end{align*}
where $A_1(n,i)$, $A_2(n,i)$ and $A_3(n,i)$ are given by
\begin{align*}
&A_1(n,i)=4(2i+4n+7)(i+n+1)^2d_i^2(n)\\[6pt]
&\qquad\qquad \quad -4(n+i)(2i+4n+11)(i+n+2)
d_{i+1}(n)d_{i-1}(n),\\[6pt]
&A_2(n,i)=(2i+4n+7)(2i+4n+5)^2d_{i+1}^2(n)\\[6pt]
&\qquad\qquad \quad
 -(2i+4n+3)(2i+4n+11)(2i+4n+7)d_i(m)d_{i+2}(n),
 \\[6pt]
&A_3(n,i)=(8i^3+40i^2+58i+32n^3+42n +80n^2+120ni+40i^2n
+64n^2i+8)\\[6pt]
&\qquad\qquad \quad  \; \cdot d_{i+1}(n)d_i(n)
-2(n+i)(2i+4n+11)(2i+4n+7)d_{i+2}(n)d_{i-1}(n).
\end{align*}
 We claim that $A_1(n,i)$, $A_2(n,i)$ and $A_3(n,i)$ are
 positive for $0 \leq i \leq n-2$.
 By the inductive hypothesis \eqref{W-assume},
 we find that for
 $0\leq i\leq n-2$,
\begin{align*}
A_1(n,i)>\,&4(2i+4n+7)(i+n+1)^2d_i^2(n)\\[6pt]
&\quad -4(n+i)(2i+4n+11)(i+n+2)\frac{(4n+2i+1)}{(4n+2i+5)} d_i^2(n)
\\[6pt]
=\,&4\frac{35+96n+72i+64ni+40n^2+28i^2}{2i+4n+5}d_i^2(n),
\end{align*}
which is positive.  From \eqref{W-assume}  it follows that for $0
\leq i \leq n-2$,
 \begin{align*}
A_2(n,i)>\,&(2i+4n+7)(2i+4n+5)^2d_{i+1}^2(n)\\[6pt]
&\quad -(2i+4n+3)(2i+4n+11)(2i+4n+7)\frac{(4n+2i+3)}{(4n+2i+7)}
d_{i+1}^2(n)\\[6pt]
=\,&(40i+80n+76)d_{i+1}^2(n),
 \end{align*}
which is also positive. By the inductive hypothesis \eqref{W-assume},
we see that  for $0\leq i \leq n-2$,
\begin{equation}\label{thm11}
d_i(n)d_{i+1}(n)>\frac{(2i+4n+5)(2i+4n+7)}{(2i+4n+3)(2i+4n+1)}
d_{i-1}(n)d_{i+2}(n).
\end{equation}
Because of \eqref{thm11}, we see that
\begin{align*}
A_3(n,i)>& (8i^3+40i^2+58i+32n^3+42n +80n^2+120ni+40i^2n +64n^2i+8)
d_{i+1}(n)d_i(n)\\[6pt]
&-2(n+i)(2i+4n+11)(2i+4n+7)
\frac{(4n+2i+3)(4n+2i+1)}{(4n+2i+5)(4n+2i+7)}d_{i+1}(n)d_i(n)\\[6pt]
 =&8\frac{5+22n+30i+44ni+24n^2+16i^2}{2i+4n+5}d_{i+1}(n)d_i(n),
\end{align*}
which is still positive for $0 \leq i \leq n-2$. Hence we deduce the
inequality \eqref{W-result-1} 
 for $0 \leq i \leq n-2$.  It remains
 to check that \eqref{W-result-1} is true for $i=n-1$,
that is,
\begin{align}\label{W-Result-2}
\frac{d_{n-1}(n+1)} {d_{n}(n+1)} <
\frac{(6n+5)d_{n}(n+1)}{(6n+9)d_{n+1}(n+1)}.
\end{align}
In view of   \eqref{Defi}, we get
\begin{align}
d_n(n+1)&=2^{-n-2}(2n+3){2n+2\choose n+1},\label{F-3}\\[6pt]
d_{n+1}(n+1)&=\frac{1}{2^{n+1}}{2n+2\choose n+1}. \label{n+1}\\[6pt]
d_n(n+2)&=\frac{(n+1)(4n^2+18n+21)}{2^{n+4} (2n+3)}{2n+4\choose
n+2}.\label{F-4}
\end{align}
Consequently,
\begin{align*}
\frac{d_{n-1}(n+1)} {d_{n}(n+1)}=
\frac{n(4n^2+10n+7)}{2(2n+1)(2n+3)}<
\frac{(2n+3)(6n+5)}{2(6n+9)}
=\frac{(6n+5)d_{n}(n+1)}{(6n+9)d_{n+1}(n+1)}.
\end{align*}
This completes the proof. \qed

 We now proceed to give a proof of (\ref{cd1}). In fact we shall prove a stronger inequality.

\begin{lem} \label{Lem-2}\label{twoline}
 Let $m\geq 2$ be a positive integer.
  For $0\leq i\leq m-1$, we have
\begin{equation}\label{tl1}
\frac{d_i(m)} {d_{i+1}(m)}>\frac{(2i+4m+5)d_i(m+1)}{(2i+4m+3)d_{i+1}(m+1)}.
\end{equation}
\end{lem}

\proof  By  Lemma \ref{strlog},
 we have for $0 \leq i \leq
m-1,$
\begin{align}\label{A-1}
d_i^2(m)>\frac{2i+4m+5}{2i+4m+1}d_{i-1}(m)d_{i+1}(m).
\end{align}
From \eqref{A-1} and  the recurrence relation  \eqref{recu1},  we
find that for $0 \leq i \leq m-1$,
\begin{align*}
& d_{i+1}(m+1)d_i(m)-\frac{2i+4m+5}{2i+4m+3}
d_{i+1}(m)d_i(m+1)\\[6pt]
&\qquad=\frac{2i+4m+5}{2(m+1)}d_{i+1}(m)d_i(m)
+\frac{i+m+1}{m+1}d_i(m)^2\\[6pt]
&\qquad\qquad -\frac{2i+4m+5}{2i+4m+3}\left(\frac{2i+4m+3}
{2(m+1)}d_i(m)d_{i+1}(m)
+\frac{i+m}{m+1}d_{i-1}(m)d_{i+1}(m)\right)\\[6pt]
&\qquad =\frac{i+m+1}{m+1}d_i^2(m)
-\frac{(4m+2i+5)(m+i)}{(4m+2i+3)(m+1)}d_{i-1}(m)d_{i+1}(m) \\[6pt]
&\qquad>\left(\frac{m+1+i}{m+1}-
\frac{(4m+2i+1)(m+i)}{(4m+2i+3)(m+1)} \right)
d_i^2(m)\\[6pt]
 &\qquad
=\frac{6m+4i+3} {(4m+2i+3)(m+1)} d_i^2(m),
\end{align*}
which is positive. This yields \eqref{tl1}, and hence the  proof is complete.
 \qed

Let us turn to the proof of (\ref{cd2}).

\noindent
  {\it Proof of (\ref{cd2}).}
  We proceed by induction on $m$. Clearly,
 the (\ref{cd2}) holds for $m=2$. We assume that it  is true
 for $n \geq 2$, that is, for $0\leq i \leq n-1$,
 \begin{align}\label{A-2}
\frac{d_i(n)}{d_{i+1}(n)}<
 \frac{d_{i+1}(n+1)}{d_{i+2}(n+1)} .
 \end{align}
It will be shown that the theorem holds for $n+1$, that is, for $0\leq i \leq n$,
 \begin{align}\label{A-3}
\frac{d_i(n+1)}{d_{i+1}(n+1)} <\frac{d_{i+1}(n+2)}{d_{i+2}(n+2)}.
 \end{align}
From the unimodality  \eqref{M}, 
 it follows that
 $d_{i}(n+1)<d_{i+1}(n+1)$ for $0
\leq i \leq  \left[\frac{n+1}{2}\right]-1$ and
$d_{i}(n+1)>d_{i+1}(n+1)$ for
 $ \left[\frac{n+1}{2}\right] \leq i
\leq n$.   From
  the recurrence  relation
\eqref{recu1}, we find that for $0 \leq i \leq
\left[\frac{n+1}{2}\right]-1$,
\begin{align*}
&d_{i+1}(n+1)d_{i+1}(n+2)-d_{i+2}(n+2)d_{i}(n+1)\\[6pt]
&\qquad=\frac{2i+4n+9}{2(n+2)}d_{i+1}^2(n+1)
+\frac{i+n+2}{n+2}d_i(n+1)d_{i+1}(n+1)\\[6pt]
&\qquad\qquad -\frac{2i+4n+11}{2(n+2)}d_i(n+1)d_{i+2}(n+1)
-\frac{i+n+3}{n+2}d_i(n+1)d_{i+1}(n+1)\\[6pt]
&\qquad=\frac{2i+4n+9}{2(n+2)}d_{i+1}^2(n+1)
-\frac{2i+4n+11}{2(n+2)}d_i(n+1)d_{i+2}(n+1)\\[6pt]
 &\qquad\qquad -\frac{1}{n+2}d_i(n+1)d_{i+1}(n+1)\\[6pt]
&\qquad>\frac{2i+4n+7}{2(n+2)}d_{i+1}^2(n+1)
-\frac{2i+4n+11}{2(n+2)}d_i(n+1)d_{i+2}(n+1),
\end{align*}
which is positive  by
 Lemma  \ref{strlog}. It follows that
for $0 \leq i \leq \left[\frac{n+1}{2}\right]-1$,
\begin{align}\label{c-1}
d_{i+1}(n+1)d_{i+1}(n+2)-d_{i+2}(n+2)d_{i}(n+1)> 0.
\end{align}
In other words, (\ref{cd2}) is valid for 
$0 \leq i \leq \left[\frac{n+1}{2}\right]-1$.

We now consider the  case $
\left[\frac{n+1}{2}\right] \leq i\leq n-1$.
 From the  recurrence relations  \eqref{recu1} and
 \eqref{recu2},  it follows that for
 $\left[\frac{n+1}{2}\right] \leq i \leq n-1$,
\begin{align*}
&d_{i+1}(n+2)d_{i+1}(n+1)-d_{i+2}(n+2)d_{i}(n+1)\\[6pt]
&\qquad=\left(\frac{(4n-2i+5)(n+i+3)}{2(n+2)(n+1-i)} d_{i+1}(n+1)
-\frac{(i+1)(i+2)}{(n+2)(n+1-i)}d_{i+2}(n+1)\right)\\[6pt]
&\qquad\qquad \times  \left(
\frac{n+1+i}{n+1}d_i(n)+\frac{4n+2i+5}{2(n+1)}d_{i+1}(n)
 \right)\\[6pt]
&\qquad\qquad -\left(\frac{n+3+i}{n+2}d_{i+1}(n+1)
+\frac{4n+2i+11}{2(n+2)}d_{i+2}(n+1) \right)\\[6pt]
&\qquad\qquad \qquad \times \left(
\frac{(4n-2i+3)(n+i+1)}{2(n+1)(n+1-i)}d_i(n)
-\frac{i(i+1)}{(n+1)(n+1-i)}d_{i+1}(n)
 \right)
\\[6pt]
&\qquad=B_1(n,i)d_{i+1}(n+1)d_i(n)+B_2(n,i)
d_{i+1}(n+1)d_{i+1}(n)\\[8pt]
&\qquad\qquad\qquad +B_3(n,i)d_{i+2}(n+1)d_i(n)+B_4(n,i)
d_{i+2}(n+1)d_{i+1}(n),
\end{align*}
where $B_1(n,i)$, $B_2(n,i)$,  $B_3(n,i)$ and $B_4(n,i)$
 are given by
\begin{align}
B_1(n,i)&=\frac{(n+i+3)(n+1+i)}{(n+2)(n+1-i)(n+1)},
\label{de-AL}\\[6pt]
B_2(n,i)&=\frac{(n+i+3)(16n^2+40n+25+4i)}{4(n+2)(n+1-i)(n+1)},
\label{de-BE}\\[6pt]
B_3(n,i)&=-\frac{(n+1+i)(41+16n^2+56n-4i)}{4(n+2)(n+1-i)(n+1)}
,\label{de-GA}\\[6pt]
B_4(n,i)&=-\frac{(i+1)(4n+5-i)}{(n+2)(n+1-i)(n+1)}. \label{de-DE}
\end{align}
Since $ \left[\frac{n+1}{2}\right] \leq i\leq n-1$, it is clear from (\ref{M}) that
$d_{i+1}(n+1) > d_{i+2}(n+1)$ and $d_{i}(n) >d_{i+1}(n)$. Thus we
get
\begin{align}
d_{i+1}(n+1)d_i(n)&>d_{i+1}(n+1)d_{i+1}(n),\label{E-1}\\[6pt]
 d_{i+1}(n+1)d_{i+1}(n)&>
d_{i+2}(n+1)d_{i+1}(n).\label{E-2}
\end{align}
Observe that $B_1(n,i)$, $B_2(n,i)$ are positive and
$B_3(n,i)$, $B_4(n,i)$  are negative.
  By
the inductive hypothesis \eqref{A-2},  \eqref{E-1} and \eqref{E-2},
 we deduce that for
 $ \left[\frac{n+1}{2}\right] \leq i\leq n-1$,
\begin{align}\label{c-2}
&d_{i+1}(n+2)d_{i+1}(n+1)-d_{i+2}(n+2)d_{i}(n+1)\nonumber\\[6pt]
&\qquad\qquad
>\left(B_1(n,i)+B_2(n,i)+
B_3(n,i)+B_4(n,i)\right)d_{i+1}(n+1)d_{i+1}(n)
\nonumber\\[6pt]
&\qquad\qquad=\frac{24n+10n^2-8ni+8i^2+13}{2(n+2)(n+1-i)(n+1)}
d_{i+1}(n+1)d_{i+1}(n)>0.
\end{align}
From the inequalities  \eqref{c-1} and \eqref{c-2}, it can be seen that
\eqref{A-3} holds for $0 \leq i \leq n-1$.
 
 We still are left with case $i=n$,  that
is,
\begin{align}\label{c-3}
\frac{d_n(n+1)}{d_{n+1}(n+1)}< \frac{d_{n+1}(n+2)}{d_{n+2}(n+2)}.
\end{align}
Applying \eqref{recu4} with $i=n+2$, we find
that
\[
\frac{d_{n}(n+1)}{d_{n+1}(n+1)}
=\frac{2n+3}{2}<\frac{2n+5}{2}=\frac{d_{n+1}(n+2)}{d_{n+2}(n+2)},
\]
as desired. This completes the proof. \qed

\section{Examples of interlacing log-concave polynomials}

Many combinatorial polynomials with only real zeros admit triangular
relations on their coefficients. The log-concavity
of  polynomials of this kind have been extensively studied. We show that several classical polynomials that are interlacing log-concave. 
To this end, we give a criterion for interlacing log-concavity based on triangular relations on
the coefficients.

\begin{thm}\label{general}
Suppose that for any $n\geq 0$,
\[
G_n(x)=\sum \limits_{k=0}^{n}T(n,k)x^k
\]
is a polynomial of degree $n$ which has only real
zeros, and suppose that the coefficients
  $T(n,k)$ satisfy  a recurrence relation of the
  following triangular form
\[
T(n,k)=f(n,k)T(n-1,k)+g(n,k)T(n-1,k-1).
\]
If
\begin{equation}\label{gen1}
\frac{(n-k)k}{(n-k+1)(k+1)}f(n+1,k+1)\leq f(n+1,k)\leq f(n+1,k+1)
\end{equation}
and
\begin{equation}\label{gen2}
g(n+1,k+1)\leq g(n+1,k)\leq \frac{(n-k+1)(k+1)}{(n-k)k}g(n+1,k+1),
\end{equation}
then the polynomials  $G_n(x)$ are interlacing
log-concave.
\end{thm}

\noindent
\proof Given the condition that  $G_n(x) $ has only real zeros, by Newton's
inequality, we have
\[
k(n-k)T(n,k)^2\geq (k+1)(n-k+1)T(n,k-1)T(n,k+1).
\]
Hence
\begin{align*}
&\quad T(n,k)T(n+1,k+1)-T(n+1,k)T(n,k+1)\\[5pt]
&\qquad =f(n+1,k+1)T(n,k)T(n,k+1)+g(n+1,k+1)T(n,k)^2\\[5pt]
&\qquad \qquad -f(n+1,k)T(n,k)T(n,k+1)-g(n+1,k)T(n,k-1)T(n,k+1)\\[5pt]
&\qquad \geq \left(f(n+1,k+1)-f(n+1,k)\right)T(n,k)T(n,k+1)\\[5pt]
&\qquad \qquad+\left(\frac{(n-k+1)(k+1)}{(n-k)k}g(n+1,k+1)
-g(n+1,k)\right)T(n,k-1)T(n,k+1),
\end{align*}
which is positive by \eqref{gen1} and \eqref{gen2}.
 It follows that
\begin{equation}\label{gen3}
\frac{T(n,k)}{T(n,k+1)}\geq \frac{T(n+1,k)}{T(n+1,k+1)}.
\end{equation}
On the other hand, we have
\begin{align*}
&\quad T(n,k+1)T(n+1,k+1)-T(n,k)T(n+1,k+2)\\[5pt]
&\qquad=f(n+1,k+1)T(n,k+1)^2+g(n+1,k+1)T(n,k)T(n,k+1)\\[5pt]
&\qquad \qquad-f(n+1,k+2)T(n,k)T(n,k+2)-g(n+1,k+2)T(n,k+1)T(n,k)\\[5pt]
&\qquad\geq \left(f(n+1,k+1)-\frac{(n-k-1)(k+1)}{(n-k)(k+2)}
f(n+1,k+2)\right)T(n,k+1)^2\\[5pt]
&\qquad \qquad+(g(n+1,k+1)-g(n+1,k+2))T(n,k+1)T(n,k).
\end{align*}
Invoking (\ref{gen1}) and (\ref{gen2}), we get
\begin{equation}\label{gen4}
\frac{T(n,k)}{T(n,k+1)}\leq
 \frac{T(n+1,k+1)}{T(n+1,k+2)}.
\end{equation}
Hence the proof is complete by combining (\ref{gen3}) and (\ref{gen4}). \qed

 Theorem \ref{general} 
 we can show that many
combinatorial polynomials which have only real zeros are interlacing
log-concave. For example, the polynomials $(x+1)^n$, $x(x+1)\cdots (x+n-1)$, the Bell polynomials, 
and the Whitney polynomials 
\[
W_{m,n}(x)=\sum \limits_{k=1}^{n}W_m(n,k)x^k,
\]
where $m$ is fixed nonnegative integer and  the coefficients $W_m(n,k)$ satisfy the recurrence relation
\[
W_m(n,k)=(1+mk)W_m(n-1,k)+W_m(n-1,k-1).
\]

To conclude, we remark that numerical evidence suggests that the Boros-Moll polynomials possess
higher order interlacing log-concavity in the spirit of the infinite-log-concavity as introduced by Moll \cite{Moll2002}.

\vspace{0.5cm}
 \noindent{\bf Acknowledgments.}  This work was supported by  the 973
Project, the PCSIRT Project of the Ministry of Education, the
Ministry of Science and Technology, and the National Science
Foundation of China.


\begin{thebibliography}{99}
\bibitem{Amdeberhan2007}
T. Amdeberhan and  V.H. Moll, A formula for a quartic integral: a
survey of old proofs and some new ones,  Ramanujan J. 18 (2008), 91--102.



\bibitem{George1999-1}
G. Boros and V.H. Moll, An integral hidden in Gradshteyn and Ryzhik,
 J. Comput. Appl. Math.  106 (1999), 361--368.

\bibitem{George1999-2}
G. Boros and V.H. Moll, A sequence of unimodal polynomials,  J.
Math. Anal. Appl.  237 (1999), 272--285.

\bibitem{George1999-3}
G. Boros and V.H. Moll, A criterion for unimodality,   Electron. J.
Combin.  6 (1999), R3.

\bibitem{George2001}
G. Boros and V.H. Moll, The double square root, Jacobi polynomials
and Ramanujan's Master Theorem, J. Comput. Appl. Math.  130 (2001),
337--344.

\bibitem{George2004}
G. Boros and V.H. Moll,  Irresistible Integrals, Cambridge
University Press, Cambridge,  2004.

\bibitem{Chen2008}
W.Y.C. Chen and E.X.W. Xia,
 The ratio monotonicity of Boros-Moll
polynomials, Math. Comp.  78 (2009), 2269--2282.

\bibitem{CG09}
W.Y.C Chen and C.C.Y. Gu, 
The reverse ultra log-concavity of the Boros-Moll polynomials, Proc. Amer. Math. Soc. 137 (2009), 3991--3998.


\bibitem{Kauers2006}
M. Kausers and P. Paule, A computer proof of Moll's log-concavity
conjecture,   Proc. Amer. Math. Soc.  135 (2007),
 3847--3856.

\bibitem{Moll2002}
V.H. Moll,  The evaluation of integrals: A personal story,  Notices
Amer. Math. Soc.  49 (2002), 311--317.

\bibitem{Moll2007}
V.H. Moll, Combinatorial sequences arising from a rational integral, Online J.
Anal. Combin. 2 (2007), \#4.


\bibitem{Wilf1992}
H.S. Wilf and D. Zeilberger, An algorithmic proof theory for
hypergeometric (ordinary and ``$q$") multisum/integral identities,
 Invent. Math. 108 (1992), 575--633.



\end{thebibliography}
\end{document}